\documentclass[reqno, a4paper]{amsart}

\usepackage{amsmath, amssymb,mathrsfs}

\usepackage{color}

\input xy
\xyoption{all}

\usepackage[colorlinks=true]{hyperref}
\hypersetup{allcolors=[rgb]{0.1,0.1,0.4}}

\newtheorem{theorem}{Theorem}[section]
\newtheorem{definition}[theorem]{Definition}
\newtheorem{proposition}[theorem]{Proposition}
\newtheorem{remark}[theorem]{Remark}
\newtheorem{lemma}[theorem]{Lemma}
\newtheorem{corollary}[theorem]{Corollary}
\newtheorem{example}[theorem]{Example}

\begin{document}

\title{Epireflective subcategories\\ and formal closure operators}

\author{Mathieu Duckerts-Antoine}
\address{Centre for Mathematics, University of Coimbra, Department of Mathematics, Apartado 3008, 3001-501 Coimbra, Portugal}
\email{mathieud@mat.uc.pt}

\author{Marino Gran}
\address{Institut de Recherche en Math\'ematique et Physique, Universit\'e catholique de Louvain, Chemin du Cyclotron, 2 bte L7.01.01, 1348 Louvain-la-Neuve, Belgium}
\email{marino.gran@uclouvain.be}

\author{Zurab Janelidze}
\address{Mathematics Division, Department of Mathematical Sciences, Stellenbosch University, Private Bag X1 Matieland, 7602, South Africa}
\email{zurab@sun.ac.za}

\thanks{The first author was partially supported by the Universit\'e catholique de Louvain, by the Centre for Mathematics of the University of Coimbra, and by the grant SFRH/BPD/98155/2013, funded by the Portuguese Government through FCT/MCTES and co-funded by the European Regional Development Fund through the Partnership Agreement PT2020. Research of the third author is supported by the South African National Research Foundation; the third author is grateful for the kind hospitality of Universit\'e catholique de Louvain.}


\maketitle

\begin{abstract} On a category $\mathscr{C}$ with a designated (well-behaved) class $\mathcal{M}$ of monomorphisms, a closure operator in the sense of D.~Dikranjan and E.~Giuli is a pointed endofunctor of $\mathcal{M}$, seen as a full subcategory of the arrow-category $\mathscr{C}^\mathbf{2}$ whose objects are morphisms from the class $\mathcal{M}$, which ``commutes'' with the codomain functor $\mathsf{cod}\colon \mathcal{M}\to \mathscr{C}$. In other words, a closure operator consists of a functor $C\colon \mathcal{M}\to\mathcal{M}$ and a natural transformation $c\colon 1_\mathcal{M}\to C$ such that $\mathsf{cod} \cdot C=C$ and $\mathsf{cod}\cdot c=1_\mathsf{cod}$. In this paper we adapt this notion to the domain functor $\mathsf{dom}\colon \mathcal{E}\to\mathscr{C}$, where $\mathcal{E}$ is a class of epimorphisms in $\mathscr{C}$, and show that such closure operators can be used to classify $\mathcal{E}$-epireflective subcategories of $\mathscr{C}$, provided $\mathcal{E}$ is closed under composition and contains isomorphisms. Specializing to the case when $\mathcal{E}$ is the class of regular epimorphisms in a regular category, we obtain known characterizations of regular-epireflective subcategories of general and various special types of regular categories, appearing in the works of the second author and his coauthors. These results show the interest in investigating further the notion of a closure operator relative to a general functor. They also point out new links between epireflective subcategories arising in algebra, the theory of fibrations, and the theory of categorical closure operators.
\end{abstract}

\section*{Introduction}

A classical result in the theory of abelian categories describes the correspondence between the localizations of a locally finitely presentable abelian category $\mathscr{C}$ and the universal closure operators on subobjects in $\mathscr{C}$ (see \cite{B1} for instance). Several related investigations in non-abelian contexts have been carried out during the last decade by several authors \cite{BCGS,BGM,BG,CDT,CG, D,DEG}. In non-abelian algebraic contexts such as groups, rings, crossed modules and topological groups, regular-epireflections are much more interesting than localizations: not only they occur more frequently but also they have strong connections with non-abelian homological algebra and commutator theory \cite{D,DEG,E}. In particular, in the pointed context of homological categories \cite{Borceux-Bourn}, the regular-epireflective subcategories were shown to bijectively correspond to a special type of closure operators on normal subobjects \cite{BG}. An analogous result was established later on in the non-pointed regular framework using closure operators on effective equivalence relations \cite{BGM}. 

By carefully examining these similar results, it appeared that the crucial idea underlying the connection between regular-epireflective subcategories and closure operators could be expressed via a suitable procedure of ``closing quotients''. Indeed, in the above mentioned situations, both normal subobjects and effective equivalence relations were ``representations'' of regular quotients. The regularity of the base category was there to guarantee the good behavior of quotients, and the additional exactness conditions only provided the faithfulness of the representation of quotients by normal subobjects/effective equivalence relations. This led to the present article where we generalize these results after introducing a general notion of a closure operator which captures both procedures --- ``closing subobjects'' and ``closing quotients''.

We now briefly describe the main content of the article. In the first section we introduce an abstract notion of a closure operator on a functor that enables us to give a common and simplified treatment of all the situations mentioned above. In the second section, we then prove our most general result, Proposition~\ref{ThA}, relating some closure operators on a specific (faithful) functor with $\mathcal{E}$-reflective subcategories, for a suitable class $\mathcal{E}$ of epimorphisms. In the third section, we make use of the concept of a form \cite{Z2,ZW1} to explain how this work extends and refines the main results concerning closure operators on normal subobjects and on effective equivalence relations. Here we also give a number of examples from algebra, where the motivation for this paper lies. For instance, specializing our results to varieties of universal algebras, we can characterize quasi-varieties in a variety via cohereditary closure operators on the ``form of congruences'', and among these varieties correspond to those closure operators which are in addition minimal. Last section is devoted to a few concluding remarks.

\section{The notion of a closure operator on a functor}\label{SecA}

\begin{definition}
A \emph{closure operator} on a functor $F\colon \mathscr{B}\to\mathscr{C}$ is an endofunctor $C\colon \mathscr{B}\to\mathscr{B}$ of $\mathscr{B}$ together with a natural transformation $c\colon 1_\mathscr{B}\to C$ such that $$FC=F\text{ and }F\cdot c=1_F.$$
A closure operator will be written as an ordered pair $(C,c)$ of the data above. A functor $C\colon \mathscr{B}\to\mathscr{B}$ alone, with the property $FC=F$, will be called an \emph{operator} on $F$. 
\end{definition}

This notion is a straightforward generalization of the notion of a categorical closure operator in the sense of D.~Dikranjan and E.~Giuli \cite{DG}. Let $\mathcal{M}$ be a class of monomorphisms in a category $\mathscr{C}$ satisfying the conditions stated in \cite{DG}. Viewing $\mathcal{M}$ as the full subcategory of the arrow-category $\mathscr{C}^\mathbf{2}$, closure operators on the codomain functor $\mathsf{cod}\colon \mathcal{M}\to\mathscr{C}$ are precisely the Dikranjan-Giuli closure operators. Almost the same is true for Dikranjan-Tholen closure operators, as defined in \cite{DT2}, which generalize Dikranjan-Giuli closure operators by simply relaxing conditions on the class $\mathcal{M}$ (see also \cite{DGT}, \cite{DT} and \cite{Tholen} for intermediate generalizations). For Dikranjan-Tholen closure operators, the class $\mathcal{M}$ is an arbitrary class of morphisms containing isomorphisms and being closed under composition with them; the closure operators are then required to satisfy an additional assumption that each component of the natural transformation $c$ is given by a morphism from the class $\mathcal{M}$ --- our definition does not capture this additional requirement.
$\mathcal{M}$ is a class of not necessarily monomorphisms already in the definition of a categorical closure operator given in \cite{DGT}; however, instead of the additional condition on a closure operator as in \cite{DT2}, there is an additional ``left-cancellation condition'' on $\mathcal{M}$ as in \cite{Tholen} (although there $\mathcal{M}$ is a class of monomorphisms) --- our definition of a closure operator for such $\mathcal{M}$ becomes the definition of a closure operator given in \cite{DGT}.

Let us remark that every pointed endofunctor $(C\colon \mathscr{B}\to\mathscr{B}, c\colon 1_\mathscr{B} \to C)$ of $\mathscr{B}$ can be viewed as a closure operator on the functor $\mathscr{B} \to \bf{1}$, where $\bf{1}$ is a single-morphism category.

In this paper we will be concerned with a different particular instance of the notion of a closure operator, where instead of a class of monomorphisms, we work with a class of epimorphisms, and instead of the codomain functor, we work with the domain functor $\mathsf{dom}\colon \mathcal{E}\to\mathscr{C}$. The motivation for the study of this type of closure operators comes from algebra, as explained in the Introduction (see the last section for some representative examples). Let us remark that these closure operators are not the same as \emph{dual closure operators} studied in \cite{DT2} (which are almost the same as \emph{coclosure operators} in the sense of \cite{DGT}). In the latter case, the functor to consider is the dual of the domain functor $\mathsf{dom}^\mathsf{op}\colon \mathcal{E}^\mathsf{op}\to\mathscr{C}^\mathsf{op}$. 

There seems to be four fundamental types of functors on which closure operators are of interest. Given a class $\mathcal{A}$ of morphisms in a category $\mathscr{C}$, regarding $\mathcal{A}$ as the full subcategory of the arrow-category of $\mathscr{C}$, these four types of functors are the domain and the codomain functors and their duals:
$$
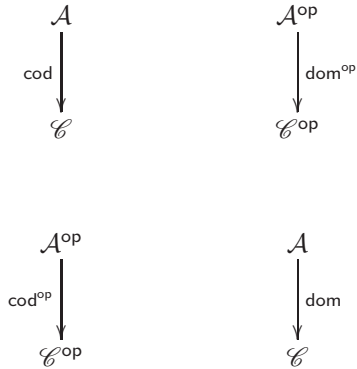
\xymatrix@=30pt{ & \mathcal{A}\ar[d]_-{\mathsf{cod}} & & \mathcal{A}^\mathsf{op}\ar[d]^-{\mathsf{dom}^\mathsf{op}} & \\  & \mathscr{C} & & \mathscr{C}^\mathsf{op} & \\ & \mathcal{A}^\mathsf{op}\ar[d]_-{\mathsf{cod}^\mathsf{op}} & & \mathcal{A}\ar[d]^-{\mathsf{dom}} & \\  & \mathscr{C}^\mathsf{op} & & \mathscr{C} & \\}$$
Horizontally, we have \emph{categorical duality}, i.e., dualizing the construction of the functor gives the other functor in the same row. Vertically, we have \emph{functorial duality}: to get the other functor in the same column, simply take the dual of the functor. The effects of closure of a morphism from the class $\mathcal{A}$ in each of the above four cases are as follows:
$$
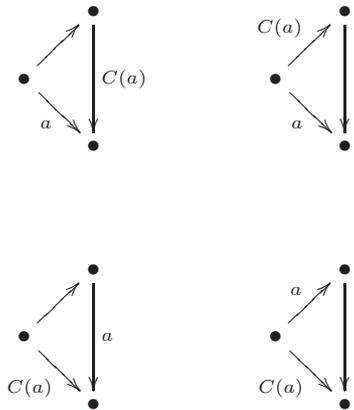
\xymatrix@=15pt{ & \bullet\ar[dd]^-{C(a)} & & & & \bullet\ar[dd] \\ \bullet\ar[ur]\ar[dr]_-{a} & & & & \bullet\ar[ur]^-{C(a)}\ar[dr]_-{a} & \\  & \bullet & & & & \bullet \\ & & & & & \\  & \bullet\ar[dd]^-{a} & & & & \bullet\ar[dd] \\ \bullet\ar[ur]\ar[dr]_-{C(a)} & & & & \bullet\ar[ur]^-{a}\ar[dr]_-{C(a)} & \\ & \bullet & & & & \bullet}$$
Note that the closure operators in the top row factorize a morphism $a$, while those in the bottom row present it as part of a factorization. This gives a principal difference between the categorical closure operators considered in the literature (which are of the kind displayed in the top row) and those that we consider in the present paper (which are of the kind displayed in the bottom row). Let us also remark that a closure operator on a poset in the classical sense can be viewed as a categorical closure operator of the bottom-right type, when we take $\mathcal{A}$ to be the class of all morphisms in the poset. Dually, the bottom-left type captures interior operators on a poset. On the other hand, categorical interior operators introduced in \cite{Vor} are not of this type (and neither it is of any of the other three types); in that case $C$ is not functorial --- it has, instead, another property that can be obtained by, yet another, ``$2$-dimensional'' duality, as remarked in Section~6 of \cite{DT2}. For a poset, the two types of closure operators in the top row become the same and they give precisely the binary closure operators in the sense of A.~Abdalla \cite{Ab}. In a poset all morphisms are both monomorphisms and epimorphisms, and it is interesting that in general, closure operators in the left column seem to be of interest when $\mathcal{A}=\mathcal{M}$ is a class of monomorphisms, and closure operators in the right column seem to be of interest when $\mathcal{A}=\mathcal{E}$ is a class of epimorphisms. In both cases the functors down to the base category are faithful. Note that another way to capture the classical notion of a closure operator on a preorder is to say that it is just a closure operator on a faithful functor $\mathscr{B} \to \bf{1}$.

Closure operators on a given functor $F$ constitute a category in the obvious way, where a morphism $n\colon (C,c)\to (C',c')$ is a natural transformation $n\colon C\to C'$ such that $n \circ c=c'$ (and consequently $F \cdot n=1_F$; note that when $F$ is a faithful functor, this last equality is equivalent to the former). We will denote this category by $\mathsf{Clo}(F)$.

For a faithful functor $F\colon \mathscr{B}\to\mathscr{C}$ from a category $\mathscr{B}$ to a category $\mathscr{C}$, an object $A$ in a fibre $F^{-1}(X)$ of $F$ will be represented by the display
$$\xymatrix@=30pt{ A\ar@{..}[d] \\ X }$$
and a morphism $A\to B$ which lifts a morphism $f\colon X\to Y$ by the display
$$\xymatrix@=30pt{ A\ar@{..}[d]\ar[r] & B\ar@{..}[d] \\ X\ar[r]_-{f} & Y }$$
Note that since the functor $F$ is faithful, it is not necessary to label the top arrow in the above display. We will also interpret this display as a statement that the morphism $f$ lifts to a morphism $A\to B$. 
When it is not clear which functor $F$ we have in mind, we will label the above square with the relevant $F$, as shown below:
$$\xymatrix@=30pt{ A\ar@{..}[d]\ar[r]\ar@{}[rd]|-{F} & B\ar@{..}[d] \\ X\ar[r]_-{f} & Y }$$

We write $A\leqslant B$ to mean 
$$\xymatrix@=30pt{ A\ar@{..}[d]\ar[r] & B\ar@{..}[d] \\ X\ar[r]_-{1_X} & X }$$
and $A\approx B$ when we also have $B\leqslant A$. In the latter case, we say that $A$ and $B$ are \emph{fibre-isomorphic}, since $A\approx B$ is equivalent to the existence of an isomorphism $A\to B$ which lifts the identity morphism $1_X$. The relation of fibre-isomorphism is an equivalence relation.

Given a faithful functor $F\colon \mathscr{B}\to\mathscr{C}$ and a morphism $f\colon X\to Y$ in $\mathscr{C}$, we will write $fA$ for the codomain of a cocartesian lifting of $f$ at $A$, when it exists. The universal property of the cocartesian lifting can be expressed as the law $$\vcenter{\xymatrix@=30pt{ A\ar@{..}[d]\ar[r] & fA\ar@{..}[d]\ar[r] & C\ar@{..}[d] \\ X\ar[r]_-{f} & Y\ar[r]_-{g} & Z }}\quad\Leftrightarrow\quad \vcenter{\xymatrix@=30pt{ A\ar@{..}[d]\ar[rr] & & C\ar@{..}[d] \\ X\ar[rr]_-{g \circ f} & & Z }}$$
More precisely, a cocartesian lifting of $f$ is the same as a lifting of $f$ satisfying the above equivalence.
Dually, we write $Cg$ for the domain of a cartesian lifting of $g$ at $C$, when it exists, and it is defined by the law $$\vcenter{\xymatrix@=30pt{ A\ar@{..}[d]\ar[r] & Cg\ar@{..}[d]\ar[r] & C\ar@{..}[d] \\ X\ar[r]_-{f} & Y\ar[r]_-{g} & Z }}\quad\Leftrightarrow\quad \vcenter{\xymatrix@=30pt{ A\ar@{..}[d]\ar[rr] & & C\ar@{..}[d] \\ X\ar[rr]_-{g \circ f} & & Z }}$$
We say $fA$ \emph{is defined} when a cocartesian lifting of $f$ at $A$ exists, and dually, we say $Cg$ \emph{is defined} when the cartesian lifting of $f$ at $C$ exists (this notation is taken from \cite{Z2,ZW1}). When $fA$ and $Cg$ are used in an equation/diagram, we interpret this equation to subsume the statement that $fA$ and $Cg$, respectively, are defined. 

Liftings of identity morphisms can be represented by vertical arrows: the display
$$\xymatrix@=30pt{ A' \\  A\ar@{..}[d]\ar[u] \\ X }$$ 
shows two objects $A$ and $A'$ in the fibre $F^{-1}(X)$, and a morphism $A\to A'$ which by $F$ is mapped to the identity morphism $1_X$.

In the case of a faithful functor $F$, the natural transformation $c$ in the definition of a closure operator is unique, when it exists, so a closure operator can be specified just by the functor $C$. In fact, it can even be given by a family $(C_X)_{X\in\mathscr{C}}$ of maps $$C_X\colon F^{-1}(X)\to F^{-1}(X),\quad A\mapsto \overline{A},$$ such that for any morphism $f\colon X\to Y$ in $\mathscr{C}$, we have the following law:
$$\vcenter{\xymatrix@=30pt{ A\ar@{..}[d]\ar[r] & B\ar@{..}[d] \\ X\ar[r]_-{f} & Y }}\quad\Rightarrow\quad \vcenter{\xymatrix@=30pt{ \overline{A}\ar[r] & \overline{B}  \\ A\ar@{..}[d]\ar[r]\ar[u] & B\ar@{..}[d]\ar[u] \\ X\ar[r]_-{f} & Y }}$$
When $F$ is faithful, $\mathsf{Clo}(F)$ is a preorder with $C\leqslant C'$ whenever $C(A)\leqslant C'(A)$ for all $A\in\mathscr{B}$. Note that the underlying pointed endofunctor of a closure operator on a faithful functor is always well-pointed, i.e., $C\cdot c =c \cdot C$. We shall say that a closure operator on a faithful functor is \emph{idempotent} when the underlying pointed endofunctor is idempotent, i.e., $C\cdot c =c \cdot C$ is an isomorphism or, equivalently, $CC \approx C$.

The notion of minimality extends directly from ordinary categorical closure operators to operators on a general faithful functor: an operator $C$ on a faithful functor $F\colon \mathscr{B}\to\mathscr{C}$ is said to be \emph{minimal} when for any two objects $A\leqslant B$ in $F^{-1}(X)$, we have that $C(B)$ is a coproduct of $B$ and $C(A)$ in the preorder $F^{-1}(X)$. Similarly as in the case of ordinary categorical closure operators, when $F^{-1}(X)$ has an initial object $0$, minimality can be reformulated by the same condition, but this time with $A=0$. Also, when $F^{-1}(X)$ has coproducts, minimality can be equivalently reformulated by requiring that $C(A+B)\approx C(A)+B$ holds for all $A,B\in F^{-1}(X)$.   

It is less trivial to extend the notion of heredity from ordinary categorical closure operators to operators on a general faithful functor. For this we will need the notion of a ``universalizer'' from \cite{Z2}, adapted to faithful functors. Consider a faithful functor $F\colon \mathscr{B}\to\mathscr{C}$ and an object $B\in F^{-1}(Y)$ 
in $\mathscr{B}$. A \emph{left universalizer} of $B$ is a morphism $f\colon X\to Y$ in $\mathscr{C}$ such that $$\xymatrix@=30pt{ A\ar@{..}[d]\ar[r] & B\ar@{..}[d] \\ X\ar[r]_-{f} & Y }$$
for any $A\in F^{-1}(X)$ and is universal with this property, i.e., for any other morphism $f'\colon X'\to Y$ such that  
$$\xymatrix@=30pt{ A'\ar@{..}[d]\ar[r] & B\ar@{..}[d] \\ X'\ar[r]_-{f'} & Y }$$
we have $f'=fu$ for a unique morphism $u$. A \emph{right universalizer} is defined dually, as a left universalizer relative to the functor $F^\mathsf{op}\colon \mathscr{B}^\mathsf{op}\to\mathscr{C}^\mathsf{op}$. 

An operator $C$ on a faithful functor $F\colon \mathscr{B}\to\mathscr{C}$ is said to be \emph{hereditary} when for any left universalizer $f\colon X\to Y$ and any object $A\in F^{-1}(X)$, when $fA$ is defined also $C_Y(fA)f$ is defined and we have $C_X(A)\approx C_Y(fA)f$.

In the case of the codomain functor $\mathsf{cod}\colon\mathcal{M}\to\mathscr{C}$, where $\mathcal{M}$ is a class of monomorphisms as in Section 3.2 of \cite{DT2}, left universalizers are simply members of the class $\mathcal{M}$ (see \cite{ZW1}) and our notion of heredity coincides with the usual one for categorical closure operators --- our heredity formula will in fact give exactly the one appearing in \cite{DT2} for $\mathcal{M}$ a pullback-stable class. On the other hand, in the case of the dual of the domain functor, $\mathsf{dom}^\mathsf{op}\colon \mathcal{E}^\mathsf{op}\to\mathscr{C}^\mathsf{op}$, for $\mathcal{E}$ having dual properties to those of $\mathcal{M}$, our notion of heredity gives the notion of coheredity of a dual closure operator in the sense of \cite{DT2}. 

Dually, an operator $C$ on a faithful functor $F\colon \mathscr{B}\to\mathscr{C}$ is said to be \emph{cohereditary} when for any right universalizer $f\colon X\to Y$ and any object $B\in F^{-1}(Y)$, when $Bf$ is defined also $f C_X(Bf)$ is defined and we have $C_Y(B)\approx f C_X(Bf)$.

\begin{remark}
The notion of \emph{weak heredity} can also be extended to an arbitrary closure operator $C$ on a faithful functor (we will, however, not make use of this notion in the present paper). Indeed, simply repeat the definition of heredity adding the assumption that $f$ is a left universalizer of $C_Y(fA)$. 
\end{remark}

\section{Closure operators for epireflective subcategories}
 
Let $\mathcal{E}$ be a class of epimorphisms in a category $\mathscr{C}$. We can view $\mathcal{E}$ as a full subcategory of the category of morphisms in $\mathscr{C}$ (the so-called ``arrow-category''), where objects are morphisms belonging to the class $\mathcal{E}$, and a morphism is a commutative square
$$\xymatrix@=30pt{ A\ar[r]^-{r} & B \\ X\ar[r]_-{f}\ar[u]^-{d} & Y\ar[u]_-{e} }$$
where $d\in\mathcal{E}$ and $e\in\mathcal{E}$ are the domain and the codomain, respectively, of the morphism. Since every morphism in the class $\mathcal{E}$ is an epimorphism, the top morphism in the above square is uniquely determined by the rest of the square. In other words, the domain functor $\mathcal{E}\to\mathscr{C}$, which maps the above square to its base, is faithful. We will use the above square to represent what we would have written as
$$\xymatrix@=30pt{ d\ar@{..}[d]\ar[r] & e\ar@{..}[d] \\ X\ar[r]_-{f} & Y }$$
for this faithful functor. 

The most standard closure operators are those that are defined on the codomain functor $\mathcal{M}\to\mathscr{C}$, where $\mathcal{M}$ is a class of monomorphisms in $\mathscr{C}$. The classical example of such a closure operator is the so-called Kuratowski closure operator on the category of topological spaces, which is given by defining the closure of an embedding $m\colon M\to X$ to be the embedding of the topological closure of the image of $m$ in $X$. In this paper we are interested in closure operators defined on the domain functor $\mathcal{E}\to\mathscr{C}$, where $\mathcal{E}$ is a class of epimorphisms. We will work with a class $\mathcal{E}$ that is closed under composition and contains identity morphisms. When $f$ is in $\mathcal{E}$, it is not difficult to see that a cartesian lifting for  
$$\xymatrix@=30pt{ ef\ar@{..}[d]\ar[r] & e\ar@{..}[d] \\ X\ar[r]_-{f} & Y }$$
under the domain functor $\mathcal{E}\to\mathscr{C}$, can be given by the square 
$$\xymatrix@=30pt{ B\ar[r]^-{1_B} & B \\ X\ar[r]_-{f}\ar[u]^-{e \circ f} & Y\ar[u]_-{e} }$$
(we could therefore write $ef = e\circ f$). We call these \emph{canonical cartesian liftings}.

\proposition\label{ThA}
Let $\mathcal{E}$ be a class of epimorphisms in a category $\mathscr{C}$ such that it contains isomorphisms and is closed under composition. There is a bijection between full $\mathcal{E}$-reflective subcategories of $\mathscr{C}$ and closure operators $C$ on the domain functor $\mathcal{E}\to\mathscr{C}$ satisfying the following conditions:
\begin{enumerate}
\item $C$ is (strictly) idempotent, i.e., for every object $e\in\mathcal{E}$ we have $C(C(e))=C(e)$ (equivalently, $CC=C$);

\item $C$ preserves canonical cartesian liftings of morphisms $f$ from the class $\mathcal{E}$, i.e., we have $$C(e)\circ f= C(e \circ f)$$ for arbitrary composable arrows $e,f\in\mathcal{E}$.
\end{enumerate}
Under this bijection, the subcategory corresponding to a closure operator consists of those objects $X$ for which $1_X=C(1_X)$, and for each object $Y$ of $\mathscr{C}$ the morphism $C(1_Y)$ gives a reflection of $Y$ in the subcategory.
\endproposition

\proof First, we show that the correspondence described at the end of the theorem gives a bijection between the objects of the poset and the preorder in question.
Let $\mathscr{X}$ be a full $\mathcal{E}$-reflective subcategory of $\mathscr{C}$, with $G$ denoting the subcategory inclusion $G\colon \mathscr{X}\to\mathscr{C}$. Consider a left adjoint $L\colon \mathscr{C}\to\mathscr{X}$ of $G$, and the unit $\eta$ of the adjunction. Since $G$ is a subcategory inclusion, each component of $\eta$ is a morphism $\eta_X\colon X\to L(X)$. Without loss of generality we may assume that the counit of the adjunction is an identity natural transformation. Then, an object $X$ of $\mathscr{C}$ belongs to the subcategory $\mathscr{X}$ if and only if $\eta_X=1_X$. We have
$$\vcenter{\xymatrix@=30pt{ A\ar@{<-}[d]_-{d}\ar[r]^-{r} & B\ar@{<-}[d]^-{e} \\ X\ar[r]_-{f} & Y }}\quad\Rightarrow\quad \vcenter{\xymatrix@=30pt{ L(A)\ar[r]^-{L(r)} & L(B)  \\ A\ar@{<-}[d]_{d}\ar[r]^-{r}\ar[u]^-{\eta_A} & B\ar@{<-}[d]^-{e}\ar[u]_-{\eta_B} \\ X\ar[r]_-{f} & Y }}$$
and this means that we can define a closure operator on the domain functor $\mathcal{E}\to\mathscr{C}$ by setting $C(e)=\eta_{\mathsf{cod}(e)} \circ e$. It is easy to see that both (a) and (b) hold for such closure operator $C$. At the same time, the full subcategory $\mathscr{X}$ of $\mathscr{C}$ can be recovered from the corresponding closure operator $C$ as the full subcategory of those objects $X$ for which $C(1_X)=1_X$. 

Given a closure operator $C$ on the domain functor $\mathcal{E}\to\mathscr{C}$, satisfying (a) and (b), we consider the full subcategory $\mathscr{X}$ of those objects $X$ in $\mathscr{C}$ such that $C(1_X)=1_X$. Consider the composite $L$ of the three functors
$$\xymatrix@=30pt{\mathcal{E}\ar[r]^-{C} & \mathcal{E}\ar[d]^-{\mathsf{cod}}\\ \mathscr{C}\ar[u]^{I}\ar[r]_-{L} & \mathscr{C} }$$ 
where $I$ maps every morphism $f\colon X\to Y$ in $\mathscr{C}$ to the morphism
$$\xymatrix@=30pt{ X\ar[r]^-{f} & Y \\ X\ar[u]^-{1_X}\ar[r]_-{f} & Y \ar[u]_-{1_Y}}$$
in the category $\mathcal{E}$, and $\mathsf{cod}$ is the codomain functor from $\mathcal{E}$ to $\mathscr{C}$. We claim that the values of $L$ lie in the subcategory $\mathscr{X}$. Indeed, we have
$$C(1_{L(X)}) \circ C(1_X)\stackrel{(b)}{=}C(1_{L(X)} \circ C(1_X))=C(C(1_X))\stackrel{(a)}{=}C(1_X)=1_{L(X)} \circ C(1_X)$$
and since $C(1_X)$ is an epimorphism, we get $C(1_{L(X)})=1_{L(X)}$. So we can consider $L$ as a functor $L\colon \mathscr{C}\to\mathscr{X}$. It follows from the construction that this functor is a right inverse of the subcategory inclusion $\mathscr{X}\to\mathscr{C}$. Since each morphism $C(1_X)\colon X\to L(X)$ is an epimorphism, it is easy to see that $L$ is a left adjoint of the subcategory inclusion $\mathscr{X}\to\mathscr{C}$, with the $C(1_X)$'s being the components of the unit of adjunction.  

To complete the proof of the bijection, it remains to show that $C(e)=C(1_{\mathsf{cod}(e)}) \circ e$. This we have by (b).
\endproof

\begin{remark}
Note that for a pullback-stable class $\mathcal{M}$ of monomorphisms in a category $\mathscr{C}$, and a closure operator on the codomain functor $\mathcal{M}$, the condition (b) (saying that cartesian liftings are preserved) can be expressed by the formula $f^{-1}(C(m))\approx C(f^{-1}(m))$, where $m,f\in\mathcal{M}$ (up to change of strict equality with isomorphism). In the special case when $m=f\circ f^{-1}(m)$ this formula expresses heredity of a closure operator (cf.~condition (HE) in \cite{DT}). A more direct link with heredity will be established further below in Lemma~\ref{LemA}.
\end{remark}

In the case of the domain functors $\mathsf{dom}\colon \mathcal{E}\to\mathscr{C}$, where objects in $\mathcal{E}$ are epimorphisms in $\mathscr{C}$, including the identity morphisms, cocartesian lifts are given by pushouts: $$\xymatrix@=30pt{ Y\ar[r] & Y+_X Z \\ X\ar[r]_-{g}\ar[u]^-{f} & Z\ar[u]_-{gf}}$$ Unlike in the case of cartesian liftings, there are in general no canonical cocartesian liftings. 

\proposition\label{Birkhoff}
Let $\mathscr{C}$ and $\mathcal{E}$ be the same as in Proposition~\ref{ThA}. If for any two morphisms $f\colon X\to Y$ and $g\colon X\to Z$ from the class $\mathcal E$, their pushout exists and the pushout injections belong to the class $\mathcal E$, then the bijection of Proposition~\ref{ThA} restricts to a bijection between:
\begin{enumerate}
\item Full $\mathcal E$-reflective subcategories $\mathscr{X}$ of $\mathscr{C}$ closed under $\mathcal E$-quotients, i.e., those having the property that for any morphism $f \colon X \to Y$ in the class $\mathcal E$ with $X$ in $\mathscr{X}$, the object $Y$ also belongs to $\mathscr{X}$.
\item Closure operators as in Proposition~\ref{ThA} having the additional property that  $$fC(e) \approx C(fe)$$
for any morphisms $f \colon X \to Y$ and $e \colon X \to E$ in the class $\mathcal E$, and moreover, when $e=C(e)$ we have $fe=C(fe)$.
\end{enumerate}
\endproposition
\proof Thanks to the bijection in Proposition~\ref{ThA}, it suffices to show that for a closure operator as in Proposition~\ref{ThA}, and the corresponding full $\mathcal{E}$-reflective subcategory $\mathscr{X}$ of $\mathscr{C}$ constructed in the proof of Proposition~\ref{ThA}, the following are equivalent:
\begin{itemize}
\item[(i)] $\mathscr{X}$ is closed in $\mathscr{C}$ under $\mathcal{E}$-quotients. 
\item[(ii)] The property on the closure operator $C$ given in (b).
\end{itemize}
Let $L$ and $\eta$ be the functor and the natural transformation that give the reflection of $\mathscr{C}$ in $\mathscr{X}$, as in the proof of Proposition~\ref{ThA}. As before, we choose $L$ and $\eta$ in such a way that an object $X$ of $\mathscr{C}$ lies in $\mathscr{X}$ if and only if $\eta_X=1_X$.

(i)$\Rightarrow$(ii):  Let $f$ and $e$ be as in (ii), and consider the morphism $g$ arising in a pushout giving a cocartesian lift of $f$ at $e$, as displayed in the bottom left square in the following diagram:
$$\xymatrix@=30pt { L(E)\ar[r] & L(E)+_X (E +_X Y)\ar[r]^-{h} & L(E +_X Y) \\ E\ar[r]^-{g}\ar[u]_-{\eta_E} & E +_X Y\ar[u]_-{g\eta_E}\ar[r]_-{1_{E +_X Y}} & E +_X Y\ar[u]_-{\eta_{E +_X Y}}\\ X\ar[r]_-{f}\ar[u]^-{e} & Y\ar[u]_-{fe} }$$
Since $\eta_{E +_X Y} \circ g=L(g) \circ \eta_E$, we get a morphism $h$ making the above diagram commute.
The top left morphism in this diagram belongs to the class $\mathcal{E}$, by the assumption on $\mathcal{E}$ given in the theorem, and so by (i), the object $L(E)+_X (E +_X Y)$ belongs to the subcategory $\mathscr{X}$. We can then use the universal property of $\eta_{E +_X Y}$ to deduce that $h$ is an isomorphism. We then get
$$fC(e)= f(\eta_E \circ e)\approx (g\eta_E)\circ (fe)\approx \eta_{E +_X Y} \circ (fe)=C(fe).$$
If $C(e)=e$, then $E$ lies in $\mathscr{X}$, and so $E+_X Y$ also lies in $\mathscr{X}$ by (i). Then $fe=C(fe)$.    

For (ii)$\Rightarrow$(i), simply take $e=1_X$ in (b).
\endproof

The next result shows how the preorder structure of closure operators is carried over to full $\mathcal{E}$-reflective subcategories, under the bijection given by Proposition~\ref{ThA}.

\proposition
Let $\mathcal{E}$ and $\mathscr{C}$ be as in Proposition~\ref{ThA}. Consider two full $\mathcal{E}$-reflective subcategories $\mathscr{X}_1$ and $\mathscr{X}_2$ of $\mathscr{C}$, and the closure operators $C_1$ and $C_2$ corresponding to them under the bijection established in Proposition~\ref{ThA}. Then $C_1\leqslant C_2$ if and only if every object in $\mathscr{X}_2$ is isomorphic to some object in $\mathscr{X}_1$.    
\endproposition   

\proof When $C_1\leqslant C_2$, for an object $X$ of $\mathscr{C}$ such that $1_X=C_2(1_X)$, we have: $$1_X\leqslant C_1(1_X)\leqslant C_2(1_X)=1_X.$$
This implies that $C_1(1_X)$ is an isomorphism, and since it is a reflection of $X$ in the subcategory $\mathscr{X}_1$, we have the morphism $C_1(1_X)$ witnessing the fact that $X$ is isomorphic to an object in $\mathscr{X}_1$. Suppose now every object in $\mathscr{X}_2$ is isomorphic to some object in $\mathscr{X}_1$. Then, for any morphism $e\colon X\to E$ from the class $\mathcal{E}$, we have $C_i(e)=C_i(1_E)\circ e$, $i\in\{1,2\}$, so to prove $C_1\leqslant C_2$, it suffices to show that $C_1(1_E)\leqslant C_2(1_E)$ for any object $E$ in $\mathscr{C}$. Since $C_2(1_E)$ is a reflection of $E$ in $\mathscr{X}_2$, its codomain lies in $\mathscr{X}_2$ and subsequently, it is isomorphic to an object lying in $\mathscr{X}_1$. Now, we can use the universal property of the reflection $C_1(1_E)$ of $E$ in $\mathscr{X}_1$ to ensure $C_1(1_E)\leqslant C_2(1_E)$.
\endproof

Let us now look at how the axioms on closure operators appearing in Propositions~\ref{ThA} and \ref{Birkhoff} are affected by isomorphism of closure operators:

\proposition\label{ProA}
Let $\mathscr{C}$ and $\mathcal{E}$ be as in Proposition~\ref{ThA}. For a closure operator $D$ on the domain functor $\mathsf{dom}\colon \mathcal{E}\to\mathscr{C}$, we have:
\begin{itemize}
\item[(i)] $D$ is isomorphic to a closure operator $C$ satisfying \ref{ThA}(a) and \ref{ThA}(b)  if and only if $DD\approx D$ and $D(e)\circ f\approx D(e\circ f)$ for arbitrary composable arrows $e,f\in\mathcal{E}$ (this last condition expresses preservation by $D$ of cartesian liftings of morphisms from the class $\mathcal{E}$). 
\end{itemize}
If further $\mathcal{E}$ satisfies the premise in Proposition~\ref{Birkhoff}, then we have: 
\begin{itemize}
\item[(ii)] $D$ is isomorphic to a closure operator $C$ satisfying the condition stated in \ref{Birkhoff}(b) if and only if $D$ satisfies the conditions stated in the second part of (i) and $D$ preserves cocartesian liftings of morphisms from the class $\mathcal{E}$, i.e., $fD(e)\approx D(fe)$ for arbitrary morphisms $f\colon X\to Y$ and $e\colon X\to E$ in the class $\mathcal{E}$.
\end{itemize}
\endproposition

\proof We first prove the only if part in each of (i) and in (ii). Suppose a closure operator $D$ is isomorphic to a closure operator $C$. If $C$ satisfies \ref{ThA}(a), then
$$D(D(e))\approx D(C(e)) \approx C(C(e))=C(e)\approx D(e)$$
for any morphism $e$ in the class $\mathcal{E}$. If $C$ satisfies \ref{ThA}(b), then  
$$D(e)\circ f\approx C(e)\circ f = C(e\circ f)\approx D(ef)$$
for arbitrary composable arrows $e,f\in\mathcal{E}$. Suppose now $\mathcal{E}$ satisfies the premise in Proposition~\ref{Birkhoff}. If $C$ satisfies the condition stated in \ref{Birkhoff}(b), then $$fD(e)\approx fC(e)\approx C(fe)\approx D(fe),$$
for arbitrary morphisms $f\colon X\to Y$ and $e\colon X\to E$ in the class $\mathcal{E}$. 

We will now prove the ``if'' parts in (i) and (ii). Consider a closure operator $D$ on the domain functor $\mathsf{dom}\colon\mathcal{E}\to\mathscr{C}$. Suppose $D$ satisfies the conditions stated in the second part of (i). Then the values of the map defined by
$$C(e)=\left\{\begin{array}{ll} e & \textrm{if }D(1_E)\textrm{ is an isomorphism,} \\ D(1_E) \circ e & \textrm{otherwise,}   \end{array}\right. $$  
are fibre-isomorphic to the values of $D$, so this gives a closure operator $C$ isomorphic to $D$. Furthermore, it is easy to see that we have
$$C(e\circ f)=C(1_E) \circ (e\circ f)  = (C(1_E) \circ e ) \circ f =C(e)\circ f,$$
as required in \ref{ThA}(b).
Since
$$D(e')\approx D(D(e'))\approx D(1_{E'})\circ D(e'),$$
for any morphism $e'\in\mathcal{E}$, where $E'$ denotes the codomain of $D(e')$, we get that
$D(1_{E'})$ is an isomorphism. We will use this fact for $e'=1_E$ in what follows. Let $e\in\mathcal{E}$ and let $E$ be the codomain of $e$. Write $E'$ for the codomain of $D(1_E)$. If $D(1_E)$ is an isomorphism, then we trivially have $C(C(e))=C(e)$. Suppose $D(1_E)$ is not an isomorphism. Since $D(1_{E'})$ is an isomorphism, we have 
$$C(C(e))=C(D(1_E)\circ e)=D(1_E) \circ e=C(e).$$
This completes the proof of the if part in (i). For the if part in (ii) we still use the same $C$. Suppose $D$ satisfies the condition stated in the second part of (ii). In view of Propositions~\ref{ThA} and \ref{Birkhoff}, it suffices to prove that for any morphism $f\colon X\to Y$ from the class $\mathcal{E}$, if $1_X=C(1_X)$ then $1_Y=C(1_Y)$. Suppose $1_X=C(1_X)$. Then $1_X\approx D(1_X)$ and since $1_Y$ is the codomain of a cocartesian lifting of $f$ at $1_X$, we have 
$$1_Y\approx fD(1_X)\approx D(1_Y),$$
which implies that $D(1_Y)$ is an isomorphism. Then $1_Y=C(1_Y)$.
\endproof

The next result shows that the formula $D(e)\circ f\approx D(e\circ f)$ which appears in \ref{ProA}(i) is in fact another way to express coheredity of a closure operator on the domain functor $\mathsf{dom}\colon\mathcal{E}\to\mathscr{C}$, and for such a closure operator, the formula $fD(e)\approx D(fe)$ in \ref{ProA}(ii) gives precisely minimality.

\lemma\label{LemA}
Let $\mathscr{C}$ and $\mathcal{E}$ be as in Proposition~\ref{ThA}. For any closure operator $C$ on the domain functor $\mathsf{dom}\colon\mathcal{E}\to\mathscr{C}$, the operator $C$ is cohereditary if and only if $C$ preserves cartesian liftings of morphisms from the class $\mathcal{E}$, i.e., the formula $C(e)\circ f\approx C(e\circ f)$ holds for all $e,f\in\mathcal{E}$. Moreover, such closure operator $C$ preserves also cocartesian liftings of morphisms from the class $\mathcal{E}$, i.e., the formula $fC(e)\approx C(fe)$ holds for all $e,f\in\mathcal{E}$ if and only if the operator $C$ is minimal.  
\endlemma

\proof
Right universalizers for $\mathsf{dom}\colon\mathcal{E}\to\mathscr{C}$ are precisely the morphisms in the class $\mathcal{E}$. So coheredity states that 
the outer rectangle in every (commutative) diagram
$$\xymatrix@=30pt{A''\ar[r] & A' \\ A\ar[r]^-{1_A}\ar[u] & A\ar[u] \\ X\ar[r]_-{f}\ar[u]_-{e\circ f}\ar@/^10pt/[uu]^-{C(e\circ f)} & Y\ar[u]^-{e}\ar@/_10pt/[uu]_-{C(e)} }$$
(with $e,f\in\mathcal{E}$) is a pushout. Since the bottom square is always a pushout, this is equivalent to the top morphism $A''\to A'$ being an isomorphism, which is equivalent to $C(e)\circ f\approx C(e\circ f)$. Now, the formula $fC(e)\approx C(fe)$ is equivalent to $fC(e)\circ f\approx C(fe)\circ f$ since $f$ is an epimorphism. When $C$ is cohereditary, it is further equivalent to $fC(e)\circ f\approx C(fe\circ f)$. The composite $fC(e)\circ f$ is in fact a coproduct of $C(e)$ and $f$ in the preorder $\mathsf{dom}^{-1}(\mathsf{dom}(f))$, while the composite $fe\circ f$ is the coproduct of $e$ and $f$ in the same preorder. Rewriting the previous formula equivalently as $C(e)+f\approx C(e+f)$ we can now recognize minimality. 
\endproof

Recall that a full subcategory $\mathscr{X}$ of a category $\mathscr{C}$ is said to be replete when it contains all objects which are isomorphic to objects already contained in $\mathscr{X}$. Recall from Section~\ref{SecA} that a closure operator $C$ is idempotent when $CC\approx C$. The work in this section leads to the following:

\begin{theorem}\label{MainTheorem}
Let $\mathcal{E}$ be a class of epimorphisms in a category $\mathscr{C}$ such that it contains isomorphisms and is closed under composition.
\begin{itemize}
\item[(a)] There is a bijection between full $\mathcal{E}$-reflective replete subcategories of $\mathscr{C}$ and isomorphism classes of cohereditary idempotent closure operators $C$ on the domain functor $\mathsf{dom}\colon \mathcal{E}\to\mathscr{C}$. 

\item[(b)] The bijection above is given by assigning to a closure operator $C$ the subcategory of $\mathscr{C}$ consisting of those objects $X$ for which $C(1_X)$ is an isomorphism, and $C(1_Y)$ gives a reflection of each object $Y$ from $\mathscr{C}$ into the subcategory.

\item[(c)] When the class $\mathcal{E}$ is closed under pushouts, the bijection above restricts to one where the subcategories are closed under $\mathcal{E}$-quotients and the closure operators are minimal.

\item[(d)] Each of the bijections above gives an equivalence between the (possibly large) poset of subcategories in question, where the poset structure is given by inclusion of subcategories, and the dual of the preorder of closure operators in question.
\end{itemize}      
\end{theorem}  

\section{Formal closure operators}

Recall that a functor is said to be amnestic when in each of its fibres, the only isomorphisms are the identity morphisms. Faithful amnestic functors were called \emph{forms} in~\cite{ZW1}. By a \emph{formal closure operator} we mean a closure operator on a form.
Any faithful functor gives rise to a form by identifying in it the fibre-isomorphic objects. The original faithful functor $F$ and the corresponding form $F'$ are related by a commutative triangle
$$\xymatrix@R=40pt{\mathscr{B}\ar[rr]^-{Q}\ar[dr]_-{F} & & \mathscr{B}'\ar[dl]^-{F'}\\ & \mathscr{C} &}$$ 
Writing $[A]_\approx$ for the equivalence class of an object $A$ in $\mathscr{B}$ under the equivalence relation of fibre-isomorphism, we have:
$$\vcenter{\xymatrix@=30pt{ A\ar[r]\ar@{..}[d]\ar@{}[dr]|-{F} & B\ar@{..}[d] \\ X\ar[r]_-{f} & Y}}\quad\Leftrightarrow\quad\vcenter{\xymatrix@=30pt{[A]_\approx\ar[r]\ar@{..}[d]\ar@{}[dr]|-{F'}  & [B]_\approx\ar@{..}[d] \\ X\ar[r]_-{f} & Y}}$$ 
The functor $Q$ is an equivalence of categories, which is surjective on objects. The above display shows what the values of $Q$ are: a morphism in $\mathscr{B}$ that fits in the left hand side display above is mapped by $Q$ to a morphism in $\mathscr{B}'$ fitting the right hand side display.
The fibres of a form are (possibly large) posets, and so the the preorder of closure operators on a form is a poset. The functor $Q$ gives rise to an equivalence of categories
$$\mathsf{Clo}(F)\approx \mathsf{Clo}(F').$$
Under this equivalence, the closure operator $C'$ on the form $F'$ associated to a closure operator $C$ on $F$ is obtained by setting $C'_X([B]_\approx)=[C_X(B)]_\approx$. Notice that since $\mathsf{Clo}(F')$ is a poset, two closure operators on $F$ correspond to the same closure operator on the associated form $F'$, under the above equivalence, if and only if they are isomorphic.

Forms associated to the domain functors $\mathsf{dom}\colon \mathcal{E}\to\mathscr{C}$ that we have been considering in this paper, were called \emph{forms of $\mathcal{E}$-quotients} in \cite{ZW1}. Theorem~\ref{MainTheorem} gives us the following:

\begin{corollary}\label{MainCorollary}
Let $\mathcal{E}$ be a class of epimorphisms in a category $\mathscr{C}$ such that it contains isomorphisms and is closed under composition. 
\begin{enumerate}
\item There is an antitone isomorphism between
the poset of full $\mathcal{E}$-reflective replete subcategories of $\mathscr{C}$ and the poset of cohereditary idempotent closure operators on the form of $\mathcal{E}$-quotients. 
\item The isomorphism above is given by assigning to a closure operator the subcategory of $\mathscr{C}$ consisting of those objects $X$ of $\mathscr{C}$ for which the initial $\mathcal{E}$-quotient is closed. 
\item When the class $\mathcal{E}$ is closed under pushouts, this isomorphism restricts to one where the subcategories are closed under $\mathcal{E}$-quotients and the closure operators are minimal.      
\end{enumerate}
\end{corollary} 

As in \cite{ZW1}, we call the form corresponding to the codomain functor $\mathcal{M}\to\mathscr{C}$, where $\mathcal{M}$ is a class of monomorphisms in a category $\mathscr{C}$, the \emph{form of $\mathcal{M}$-subobjects}. A normal category in the sense of~\cite{Z} is a regular category~\cite{BGV} which is pointed and in which every regular epimorphism is a normal epimorphism. In a normal category, for the class $\mathcal{E}$ of normal epimorphisms and the class $\mathcal{M}$ of normal monomorphisms, the form of $\mathcal{E}$-quotients is isomorphic to the form of $\mathcal{M}$-subobjects, via the usual kernel-cokernel correspondence between normal quotients and normal subobjects. Corollary~\ref{MainCorollary} then gives:

\begin{theorem}\label{ThB}
There is an antitone isomorphism between the poset of full normal-epi-reflective replete subcategories of a normal category $\mathscr{C}$ and the poset of cohereditary idempotent closure operators on the form of normal subobjects. It is given by assigning to a closure operator the subcategory of $\mathscr{C}$ consisting of those objects $X$ of $\mathscr{C}$ for which the null subobject of $X$ is closed. Furthermore, when pushouts of normal epimorphisms along normal epimorphisms exist, this isomorphism restricts to one where the subcategories are closed under normal quotients and the closure operators are minimal.     
\end{theorem}

This recovers Theorem 2.4 and Proposition 3.4 from \cite{BG}, and moreover, slightly generalizes and refines them. Let us explain this in more detail. First of all, we remark that an \emph{idempotent closure operator on kernels} defined in \cite{BG} is the same as an idempotent closure operator in the sense of the present paper, on the form of normal subobjects. The context in which these closure operators are considered in~\cite{BG} is that of a homological category 
~\cite{Borceux-Bourn}, which is the same as a pointed regular protomodular category \cite{Bou91}. Theorem 2.4 in \cite{BG} establishes, for a homological category, a bijection between the so-called homological closure operators and normal-epi-reflective subcategories (which in \cite{BG} are simply called epi-reflective subcategories). This bijection is precisely the one established by the first half of Theorem~\ref{ThB} above (so, homological = cohereditary + idempotent). As this theorem shows, the bijection is there more generally for any normal category (a homological category is in particular a normal category, but the converse is not true). 

\begin{example}
\begin{enumerate}
\item 
Let $\mathsf{CRng}$ be the category of commutative (not necessarily unital) rings, and let $\mathsf{RedCRng}$ be its full reflective subcategory of reduced commutative rings (i.e., with no non-zero nilpotent element). The category $\mathsf{CRng}$ is homological and the homological closure operator associated with the corresponding reflection can be described explicitly, and it actually gives the well known notion of \emph{nilradical} of an ideal. Indeed, for any ideal $I$ of a commutative ring $A$, its closure in $A$ is its nilradical (see~\cite{DEG})
\[
  \sqrt{I}= \{ a \in A \mid \exists_{n \in \mathbb{N}}\; a^n \in I\}.
 \]
\item Consider the category $\mathsf{Grp}(\mathsf{Top})$ of topological groups and its full reflective subcategory $\mathsf{Grp}(\mathsf{Haus})$ of Hausdorff groups. The category $\mathsf{Grp}(\mathsf{Top})$ is homological and, under the closure operator corresponding to the reflective subcategory $\mathsf{Grp}(\mathsf{Haus})$, the closure of a normal subgroup $H$ of a topological group $A$ is simply given by its topological closure $\overline{H}$ in $A$ (see \cite{BG}).
\end{enumerate}
\end{example}

The last part of Theorem~\ref{ThB} similarly captures Proposition~3.4 from \cite{BG} characterizing Birkhoff subcategories \cite{JK} of a semi-abelian category. Once again, it reveals a more general context where the result can be stated, and namely that of a normal category with pushouts of normal epimorphisms along normal epimorphisms in the place of a semi-abelian category \cite{JMT}. Thus, in particular, the characterization remains valid in any ideal determined category \cite{JMTU}.

\begin{example}
\begin{enumerate}
\item Consider the category $\mathsf{Grp}(\mathsf{HComp})$ of compact Hausdorff groups and its full reflective subcategory $\mathsf{Grp}(\mathsf{Prof})$ of profinite groups. Here, the closure of a normal subgroup $H$ of a compact Hausdorff group $A$ is precisely the group-theoretic product $H \cdot \Gamma_A(1)$, where $\Gamma_A(1)$ is the connected component in $A$ of the neutral element $1$ (see \cite{BG}).
\item Let $\mathsf{PXMod}$ be the category of precrossed modules and $\mathsf{XMod}$ its full reflective subcategory of crossed modules. We recall that a \emph{precrossed module} is a group homomorphism $\alpha \colon A \to B$ together with an action of the group $B$ on $A$, denoted by ${}^ba$,  such that $\alpha({}^ba)=b\cdot \alpha(a)\cdot b^{-1}$ for all  $a\in A$ and  $b \in B$. A \emph{crossed module} is a precrossed module such that ${}^{\alpha(a_1)}a_2=a_1\cdot a_2 \cdot a_1^{-1}$ for all $a_1,a_2 \in A$. A morphism $f$ of precrossed modules from $\mu \colon M \to N$ to $\alpha \colon A\to B$  is an equivariant pair $(f_1,f_0)$ of group homomorphisms making the diagram
\[
\xymatrix@=30pt{
M \ar[r]^-{f_1} \ar[d]_-{\mu}& A\ar[d]^-{\alpha} \\
N \ar[r]_-{f_0}& B
}
\]
commute.
The category $\mathsf{XMod}$ is a Birkhoff subcategory of the semi-abelian category $\mathsf{PXMod}$ (see~\cite{E} for more details). Given a normal sub-precrossed module $\mu\colon M \to N $ of $\alpha\colon A\to B$, its closure is given by the supremum $\mu \vee \langle A,A \rangle$ of $\mu$ and  $\langle A,A \rangle$ considered as normal sub-precrossed modules of $\alpha$, where $\langle A,A \rangle$ arises as the normal subgroup of $A$ generated by 
\[
\{ {}^{\alpha(a_1)}a_2\cdot a_1\cdot a_2^{-1}\cdot a_1^{-1} \mid a_1,a_2 \in A\}.
\] 
\end{enumerate}
\end{example}

For a category $\mathscr{C}$, consider the full subcategory $\mathscr{B}$ of the category of parallel pairs of morphisms in $\mathscr{C}$, consisting with those parallel pairs of morphisms which arise as kernel pairs of a morphism $f$ (i.e., projections in a pullback of $f$ with itself). Thus, a morphism in $\mathscr{B}$ is a diagram
 \[
 \xymatrix@=30pt{
 R \ar[r]^-{g} \ar@<0.5ex>[d]^-{r_2} \ar@<-0.5ex>[d]_-{r_1}  & S\ar@<0.5ex>[d]^-{s_2} \ar@<-0.5ex>[d]_-{s_1} \\
 X \ar[r]_-{f}  & Y
 }
 \]
where $(R,r_1,r_2)$ and $(S,s_1,s_2)$ are kernel pairs, and we have $$f \circ r_1=s_1\circ g\text{ and }f \circ r_2=s_2\circ g.$$ Assigning to the above diagram the base morphism $f$ defines a (faithful) functor $\mathscr{B}\to\mathscr{C}$. The form corresponding to the functor will be called the \emph{congruence form} of $\mathscr{C}$ (when $\mathscr{C}$ is a variety of universal algebras, its fibres are isomorphic to congruence lattices of algebras). For a regular category, the congruence form is isomorphic to the form of regular quotients, and Corollary~\ref{MainCorollary} can be rephrased as follows:

\begin{theorem}\label{ThC}
There is an antitone isomorphism between the poset of full regular-epi-reflective replete subcategories of a regular category $\mathscr{C}$ and the poset of cohereditary idempotent closure operators on the congruence form. It is given by assigning to a closure operator the subcategory of $\mathscr{C}$ consisting of those objects $X$ of $\mathscr{C}$ for which the smallest congruence on $X$ is closed. Furthermore, when pushouts of regular epimorphisms along regular epimorphisms exist, this isomorphism restricts to one where the subcategories are closed under regular quotients and the closure operators are minimal.     
\end{theorem}

The first part of the theorem above recovers Theorem~2.3 from~\cite{BGM}. Idempotent closure operators on the congruence form of a regular category are the same as idempotent closure operators on effective equivalence relations in the sense of~\cite{BGM}. The condition of coheredity defines precisely the effective closure operators in the sense of~\cite{BGM}. The last part of the above theorem includes Proposition 3.6 from~\cite{BGM} as a particular case. 

\begin{example} Consider the category $\mathsf{Qnd}$ of quandles. Recall that a \emph{quandle} is a set $A$ equipped with two binary operations $\lhd$ and $\lhd^{-1}$ such that the following identities hold, for all $a,b,\in A$:
\[
a \lhd a = a,\quad (a \lhd b) \lhd^{-1} b = a = (a \lhd^{-1} b) \lhd b,\quad (a \lhd b) \lhd c = (a \lhd c) \lhd (b \lhd c).
\]
A \emph{quandle homomorphism} is a function preserving both operations. A quandle is \emph{trivial} when $a\lhd b = a$ and $a\lhd^{-1} b=a$ for all $a,b \in A$. The category $\mathsf{Qnd}$ is regular and the full subcategory $\mathsf{Qnd}^\star$ of trivial quandles is a Birkhoff subcategory  of $\mathsf{Qnd}$. In the category $\mathsf{Qnd}$, an effective equivalence relation $R$ on a quandle $A$ is a \emph{congruence} on $A$, namely an equivalence relation on the underlying set of $A$ which is compatible with the quandle operations of $A$. Given two elements $a,b\in A$, we write $a \sim_A b$ if there exist a chain of elements $a_1,\hdots, a_n \in A$ such that $(\hdots(a \lhd^{\alpha_1} a_1) \lhd^{\alpha_2} \dots ) \lhd^{\alpha_n} a_n = b$ where $\alpha_i\in\{-1,1\}$ for all $1\leqslant i \leqslant n$.  This defines a congruence $\sim_A$ on $A$. Given a congruence $R$ on $A$, its closure (relative to the reflective subcategory $\mathsf{Qnd}^\star$ of $\mathsf{Qnd}$) is given by the composite of congruences 
\[
R \circ {\sim_A} = \{ (a,b) \in A\times A \mid \exists_{c \in A}\;(a \sim_A c\; \wedge \;  c \mathrel{R} b)\}.
\] 
For more details, the reader is referred to \cite{EG}.
\end{example}

Finally, let us remark that Theorem~\ref{ThB} can be deduced already from Theorem~\ref{ThC}, since for a normal category the form of normal subobjects is isomorphic to the congruence form.

Applying Theorem~\ref{ThC} in the case when $\mathscr{C}$ is a variety of universal algebras, the first part of the theorem gives a characterization of quasi-varieties of algebras in the variety, and the second part --- subvarieties of the variety.

\section{Concluding remarks}

The notion of a categorical closure operator has a long history. Its origins lie in classical category theory, where they appear as universal closure operators (see e.g.~\cite{B1}) arising in the study of abelian categories and topoi. The notion introduced in \cite{DG} led to establishing the study of categorical closure operators as a separate subject. This development eventually inspired a new way of thinking: take a structural presentation of a topological space and turn it into a structure on a category. In particular, ``interior operators'' were introduced in \cite{Vor} and ``neighborhood operators'' were introduced in \cite{HolSla}. During the last few years the third author has proposed in a number of his talks that it may be worthwhile to define and study these structures relative not only to a category equipped with a class $\mathcal{M}$ of monomorphisms, but relative to a category equipped with a more general structure, such as a cover relation in the sense of \cite{J08,J09} or a form in the sense of \cite{Z2,ZW1} (which can be seen as a generalization of a cover relation). 

In the present paper we have tried to illustrate worthiness of studying closure operators relative to a form. In particular, we showed that it opens a way to a new type of closure operators, and namely, closure operators on the forms of $\mathcal{E}$-quotients. As we have seen, such closure operators capture epireflective subcategories through the notions of idempotence, coheredity, and minimality; this gives new and a more general perspective on the work carried out in \cite{BG} and \cite{BGM}. It would be interesting to find a similar application-based motivation for extending, to the context of forms, interior and neighborhood operators, as well as their generalizations, and also to explore usefulness of formal closure operators further. As we have seen in this paper, all standard properties of categorical closure operators generalize to these closure operators.

\subsection*{Acknowledgement.}
We would like to thank the anonymous referee for the useful report on the first version of this paper.

\end{document}